\pdfoutput=1

\documentclass[graybox]{svmult}

\usepackage{type1cm}	    % activate if the above 3 fonts are
\usepackage{makeidx}	     % allows index generation
\usepackage{amscd}
\usepackage[active]{srcltx}
\usepackage{verbatim}
\usepackage{epsfig,graphicx,color}
\usepackage{graphicx}
\usepackage{mathtools,xparse}
\usepackage[english]{babel}
\usepackage{upgreek}

\usepackage[all]{xy}

\usepackage{dsfont}
\usepackage{multicol}	     % used for the two-column index
\usepackage[bottom]{footmisc}% places footnotes at page bottom

\usepackage{pgf, tikz, tikz-cd}
\usepackage{dsfont}
\usetikzlibrary{matrix}
\usetikzlibrary{arrows}
\usetikzlibrary{patterns}
\usepackage{braids}

\usepackage{capt-of}
\usepackage{float}

\usepackage{epigraph}
\usepackage[unicode,bookmarks, pdftex, backref = page]{hyperref}
\definecolor{newblue}{RGB}{0,102,204}
\hypersetup{colorlinks=true,citecolor=newblue,linkcolor=newblue,
	urlcolor=magenta,
pdfpagemode=UseNone, breaklinks=true}

\usepackage{mathtools}
\mathtoolsset{showonlyrefs=true}
\usepackage{amscd}
\usepackage[active]{srcltx}
\usepackage{verbatim}
\usepackage{epsfig,color}
\usepackage{graphicx}
\usepackage{mathtools,xparse}

\usepackage[english]{babel}

\usepackage{upgreek}
\usepackage{dsfont}

\renewcommand{\to}{\longrightarrow}
\newcommand{\nt}{\mathsf{t}}
\newcommand{\nr}{\mathsf{r}}
\newcommand{\z}{\mathsf{z}}
\newcommand{\x}{\mathsf{x}}

\renewcommand{\S}{\mathfrak{S}}

\usepackage{amsmath,amsfonts, mathrsfs, mathtools,  fancyhdr}
\usepackage{amssymb}

\makeindex	       % used for the subject index
\usepackage{newtxtext}	 % selects Times Roman as basic font
\usepackage{newtxmath}

\begin{document}

\title*{Finite quotients of surface braid groups and double Kodaira
fibrations}
\author{Francesco Polizzi and Pietro Sabatino}
\institute{Francesco Polizzi \at Dipartimento di Matematica e Informatica,
	Universit\`{a} della Calabria, Ponte Bucci Cubo 30B, 87036
	Arcavacata di
	Rende (Cosenza), Italy, \email{francesco.polizzi@unical.it}
	\and Pietro Sabatino \at Via Val Sillaro 5, 00141 Roma
\email{pietrsabat@gmail.com}}

\maketitle

\epigraph{\itshape To Professor Ciro Ciliberto on the occasion of his
70th birthday}

\abstract{Let $\Sigma_b$ be a closed 
	Riemann surface of genus $b$. We give an account of some results obtained in the recent papers
	\cite{CaPol19, Pol20, PolSab21}
	and concerning what we call here \emph{pure braid quotients},
	namely non-abelian finite groups appearing as quotients of the pure braid group
	on two strands $\mathsf{P}_2(\Sigma_b)$. We also explain how these groups 
	can be used
	in order to provide new constructions of double Kodaira fibrations.
}

\section{Introduction} \label{sec:intro}

A \emph{Kodaira fibration} is a smooth, connected holomorphic fibration
$f_1 \colon S \to B_1$, where $S$ is a compact complex surface and $B_1$
is a compact complex curve, which is not isotrivial (this means that not
all its fibres are biholomorphic to each others). The genus $b_1:=g(B_1)$ is
called the \emph{base genus} of the fibration, whereas the genus $g:=g(F)$,
where $F$ is any fibre, is called the \emph{fibre genus}. If a surface $S$
is the total space of a Kodaira fibration, we will call it a \emph{Kodaira
fibred surface}.  For every Kodaira fibration we have $b_1 \ge 2$ and $g \geq 3$,
see \cite[Theorem 1.1]{Kas68}. Since the fibration is smooth, the condition
on the base genus implies that $S$ contains no rational or elliptic curves;
hence it is minimal and, by the sub-additivity of the Kodaira dimension,
it is of general type, hence algebraic.

Examples of Kodaira fibrations were originally constructed in \cite{Kod67,
At69} in order to show that, unlike the topological Euler characteristic,
the signature $\sigma$ of a real manifold is not multiplicative for
fibre bundles. In fact, every Kodaira fibred surface $S$ satisfies
$\sigma(S) >0$, see for example the introduction of \cite{LLR17}, whereas
$\sigma(B_1)=\sigma(F)=0$, and so $\sigma(S) \neq \sigma(B_1)\sigma(F)$.

A \emph{double Kodaira surface} is a compact complex surface $S$, endowed
with a \emph{double Kodaira fibration}, namely a surjective, holomorphic map
$f \colon S \to B_1 \times B_2$ yielding, by composition with the natural
projections, two Kodaira fibrations $f_i \colon S \to B_i$, $i=1, \,2$.

The purpose of this article is to give an account of recent results,
obtained in
the series of papers \cite{CaPol19, Pol20, PolSab21}, concerning the
construction
of some double Kodaira fibrations (that we call \emph{diagonal}) by means
of group-theoretical methods. Let us start by
introducing the needed terminology. Let $b \geq 2$ and  $n \geq 2$ 
be two positive integers, and let $\mathsf{P}_2(\Sigma_b)$ be the
pure braid group on two strands on a closed Riemann surface of genus $b$. 
We say that a finite group $G$ is a \emph{pure braid quotient}
of type $(b, \, n)$ if there exists a group epimorphism
\begin{equation} \label{eq:intro-pure-braid-quotient}
	\varphi \colon \mathsf{P}_2(\Sigma_b) \to G
\end{equation}
such that $\varphi(A_{12})$ has order $n$, where $A_{12}$ is the
braid corresponding, via the isomorphism $\mathsf{P}_2(\Sigma_b) \simeq
\pi_1(\Sigma_b \times \Sigma_b - \Delta)$, to the homotopy class in $\Sigma_b \times \Sigma_b - \Delta$ of a loop in
$\Sigma_b \times \Sigma_b$  ``winding once" around the diagonal $\Delta$. 
Since $A_{12}$ is a commutator in $\mathsf{P}_2(\Sigma_b)$ and $n \geq 2$, it follows that
every pure braid quotient is a non-abelian group, see Remark \ref{rmk:purebraidquotient_nonabelian}.

By Grauert-Remmert's extension theorem together with Serre's GAGA, the existence
of a pure braid quotient  as in \eqref{eq:intro-pure-braid-quotient} is
equivalent to the existence of a Galois cover $\mathbf{f} \colon S \to
\Sigma_b \times \Sigma_b$, branched over $\Delta$ with branching
order $n$. After Stein factorization, this yields in turn a diagonal double
Kodaira fibration $f \colon S \to \Sigma_{b_1} \times \Sigma_{b_2}$. We have
$\mathbf{f}=f$, i.e. no Stein factorization is needed, if and only if $G$
is a \emph{strong} pure braid quotient, an additional condition explained
in Definition \ref{def:strong}.

We are now in a position to state our first results, see Theorems \ref{thm:1},
\ref{thm:2}, \ref{thm:3}:
\begin{itemize}
	\item If $b \geq 2$ is an integer and $p \geq 5$ is a prime number,
		then
		both extra-special $p$-groups of order $p^{4b+1}$ are
		non-strong pure braid
		quotients of type $(b, \, p)$.
	\item If $b \geq 2$ is an integer and $p$ is a prime number
		dividing $b+1$,
		then both extra-special $p$-groups of order $p^{2b+1}$
		are pure braid quotients of type $(b, \, p)$.
	\item If a finite group $G$ is a pure braid quotient, then $|G|
		\geq 32$, with equality holding if and only if $G$ is
		extra-special. Moreover, in
		the last case, we can explicitly compute the number of
		distinct quotients
		maps of type \eqref{eq:intro-pure-braid-quotient}, up to
		the natural action
		of $\operatorname{Aut}(G)$.
\end{itemize}
We believe that such results are significant because, although we know
that $\mathsf{P}_2(\Sigma_b)$ is residually $p$-finite for all $p \geq 2$
(see \cite[pp. 1481-1490]{BarBel09}), it is usually tricky to explicitly
describe its non-abelian finite quotients.

The geometrical counterparts of  the above
group-theoretical statements allow us to construct infinite
families of double Kodaira fibrations with interesting numerical properties,
for instance having slope greater than $2 + 1/3$ or signature equal to $16$,
see Theorems \ref{thm:slope}, \ref{thm:lim-sup}, \ref{thm:3-geometric}:
\begin{itemize}
	\item Let $f \colon S_p  \to \Sigma_{b'} \times \Sigma_{b'}$ be
		the diagonal
		double Kodaira fibration associated with a non-strong pure
		braid quotient
		$\varphi \colon \mathsf{P}_2(\Sigma_2) \to G$ of type  $(2,
		\, p)$, where
		$G$ is an extra-special $p$-group $G$ of order $p^9$ and $b'=p^4+1$. Then
		the maximum
		slope $\nu(S_p)$ is attained for precisely two values of
		$p$, namely
		\begin{equation}
			\nu(S_{5}) = \nu(S_{7})= 2 + \frac{12}{35}.
		\end{equation}
		Furthermore, $\nu(S_{p}) > 2 + 1/3$ for all $p \geq 5$. More
		precisely,
		if $p \geq 7$ the function $\nu(S_{p})$ is strictly
		decreasing and
		\begin{equation*}
			\lim_{p \rightarrow +\infty} \nu(S_{p}) = 2 +
			\frac{1}{3}.
		\end{equation*}
	\item Let $\Sigma_b$ be any closed Riemann surface of genus $b$. Then
		there exists a double Kodaira fibration $f \colon S \to
		\Sigma_b \times
		\Sigma_b$. Moreover, denoting by $\kappa(b)$ the number of
		such fibrations,
		we have
		\begin{equation*}
			\kappa(b) \geq \boldsymbol{\upomega}(b+1),
		\end{equation*}
		where $\boldsymbol{\upomega} \colon \mathbb{N} \to \mathbb{N}$
		stands for the
		arithmetic function counting the number of distinct prime
		factors of a
		positive integer. In particular,
		\begin{equation*}
			\limsup_{b \rightarrow + \infty} \kappa(b) = + \infty.
		\end{equation*}
	\item Let $G$ be a finite group and $\mathbf{f} \colon S \to
		\Sigma_b \times
		\Sigma_b$ be a Galois cover, with Galois group $G$, branched
		over the diagonal
		$\Delta$ with branching order $n$. Then $|G|\geq 32$,
		and equality  holds
		if and only if $G$ is extra-special. If $G$ is extra-special
		of order $32$
		and $(b, \, n)=(2, \, 2)$, then $\mathbf{f} \colon S \to
		\Sigma_2 \times \Sigma_2$ is a diagonal double Kodaira fibration 
		such that
		\begin{equation}
			b_1 = b_2=2, \quad g_1=g_2=41, \quad
			\sigma(S)=16.
		\end{equation}
\end{itemize}
As a consequence of the last result, we obtain a sharp lower bound
for the signature of a diagonal double Kodaira fibration, see Theorem
\ref{thm:signature}:
\begin{itemize}
	\item Let $f \colon S \to \Sigma_{b_1} \times \Sigma_{b_2}$ be
		a diagonal
		double
		Kodaira fibration, associated with a pure braid quotient
		$\varphi
		\colon \mathsf{P}_2(\Sigma_b) \to G$ of type $(b, \, n)$. Then
		$\sigma(S) \geq 16$, and equality
		holds precisely when $(b, \, n)=(2, \, 2)$ and $G$ is
		an extra-special
		group of order $32$.
\end{itemize}
Note that our methods show that \emph{every} curve of genus $b$ (and
not only some special curve with extra automorphisms) is the basis of a
(double) Kodaira fibration and that, in addition, the number of distinct
Kodaira fibrations over a fixed base can be arbitrarily large. Furthermore,
\emph{every} curve of genus $2$ is the base of a (double) Kodaira fibration
with signature $16$ and this provides, to our knowledge, the first example
of positive-dimensional family of (double) Kodaira fibrations with small
signature.

The aforementioned examples with signature $16$ also provide new ``double solutions'' 
to a problem, posed by
G. Mess and included  in Kirby's problem list in low-dimensional topology, see
\cite[Problem 2.18 A]{Kir97}, asking what is the smallest number $b$ for
which there exists a real surface bundle over a real surface with base
genus $b$ and non-zero signature. We actually have $b=2$, also for double
Kodaira fibrations, see Theorem \ref{thm:real-manifold}:

\begin{itemize}
	\item Let  $S$ be 
		double
		Kodaira surface, associated with a pure braid quotient
		$\varphi \colon
		\mathsf{P}_2(\Sigma_b) \to G$ of type $(2, \, 2)$, where $G$
		is  an
		extra-special group of order $32$. Then the real
		manifold $X$
		underlying
		$S$ is a closed, orientable $4$-manifold of signature
		$16$ that
		can be realized as a real surface bundle over a real
		surface of
		genus $2$,
		with fibre  genus $41$, in two different ways.
\end{itemize}

In fact, it is an interesting question whether $16$ and $41$ are the minimum
possible values for the signature and the fibre genus of a (not necessarily
diagonal) double Kodaira surface $f \colon S \to \Sigma_2 \times \Sigma_2$,
but we will not develop this point here.

The above results paint a rather clear picture regarding pure braid
quotients and the relative diagonal double Kodaira fibrations when \( | G |
\le 32 \). It is natural then to investigate further this topic for \( | G
| > 32 \), and indeed this paper also contains the following new result, 
see  Theorem \ref{thm:new}:
\begin{itemize}
	\item If $G$ is a finite group with $32 < |G| <64$, 
		then $G$ is not a pure braid quotient.
\end{itemize}
We provide only a sketch of the proof, which is based on calculations
performed by means of the computer algebra system \verb|GAP4|, see \cite{GAP4}; 
the details will appear in a forthcoming paper.
\medskip

\textbf{Acknowledgments.} F. Polizzi was partially supported by
GNSAGA-INdAM. Both authors thank A. Causin for drawing the figures.

\medskip
\indent \textbf{Notation and conventions.}
The order of a finite group $G$ is denoted by $|G|$. If $x \in G$, the
order of $x$ is denoted by $o(x)$. The subgroup generated by $x_1, \ldots,
x_n \in G$ is denoted by $\langle x_1, \ldots, x_n \rangle$.
The center of $G$ is denoted by $Z(G)$ and the centralizer of an element
$x \in G$ by $C_G(x)$. If $x, \, y \in G$,
their commutator is defined as $[x,y]=xyx^{-1}y^{-1}$.
We denote both the cyclic group of order $p$ and the field with $p$
elements  by $\mathbb{Z}_p$. We use sometimes the \verb|IdSmallGroup(G)|
label from \verb|GAP4| list of small groups. For instance,
$\mathsf{S}_4=G(24, \, 12)$ means that $\mathsf{S}_4$ is the
twelfth group of order $24$ in this list.

\section{Pure surface braid groups and finite braid quotients}
\label{sec:pure-braids}

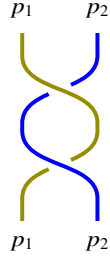
\begin{figure}[H]
\centering
\begin{tikzpicture}
\braid[number of strands=2, style strands={1}{olive}, style strands={2}{blue}, line width=1.5pt]
(braid) a_1 a_1;
\node at (braid-1-s)[yshift = 0.3cm]{$p_1$};
\node at (braid-2-s)[yshift = 0.3cm]{$p_2$};
\node at (braid-1-e)[yshift = -0.3cm]{$p_1$};
\node at (braid-2-e)[yshift = -0.3cm]{$p_2$};
\end{tikzpicture}
\medskip
\caption{A pure braid on two strands} \label{fig:braid}
\end{figure}

Let $\Sigma_b$ be a closed Riemann surface of genus $b \geq 2$, and let
$\mathscr{P} = (p_1, \, p_2)$ be an ordered pair of distinct points on
$\Sigma_b$. A \emph{pure geometric braid} on $\Sigma_b$ based at $\mathscr{P}$
is a pair $(\alpha_1, \, \alpha_2)$ of paths $\alpha_i \colon [0, \, 1]
\to \Sigma_b$ such that
\begin{itemize}
	\item $\alpha_i(0) = \alpha_i(1)=p_i \quad \textrm{for all }i \in
		\{1,\, 2\}$
	\item the points $\alpha_1(t), \, \alpha_2(t) \in \Sigma_b$
		are pairwise
		distinct for all $t \in [0, \, 1],$
\end{itemize}
see Figure \ref{fig:braid}.

\begin{definition} \label{def:braid}
	The \emph{pure braid group} on two strands on $\Sigma_b$ is the group
	$\mathsf{P}_{2}(\Sigma_b)$ whose elements are the pure braids
	based at
	$\mathscr{P}$ and whose operation is the usual concatenation of paths,
	up to homotopies among braids.
\end{definition}
It can be shown that $\mathsf{P}_{2}(\Sigma_b)$ does not depend on the choice
of the set $\mathscr{P}=(p_1, \, p_2)$, and that there is an isomorphism
\begin{equation} \label{eq:iso-braids}
	\mathsf{P}_{2}(\Sigma_b) \simeq \pi_1(\Sigma_b \times \Sigma_b
		- \Delta,
	\, \mathscr{P})
\end{equation}
where $\Delta \subset \Sigma_b \times \Sigma_b$ is the diagonal.

The group $\mathsf{P}_{2}(\Sigma_b)$  is finitely presented for all $b$,
and explicit presentations can be found in \cite{Bel04, Bir69, GG04,
S70}. Here we follow the approach in \cite[Sections 1-3]{GG04}, referring
the reader to that paper for further details.

\begin{proposition}[{\cite[Theorem 1]{GG04}}] \label{prop:split-braid}
	Let $p_1, \, p_2 \in \Sigma_b$, with $b \geq 2$. Then the map
	of pointed
	topological spaces given by the projection onto the first component
	\begin{equation} \label{eq:proj-first}
		(\Sigma_b \times \Sigma_b - \Delta, \, \mathscr{P}) \to
		(\Sigma_b, \, p_1)
	\end{equation}
	induces a split short exact sequence of groups
	\begin{equation} \label{eq:split-braid}
		1 \longrightarrow \pi_1(\Sigma_b - \{p_1\}, \, p_2)
		\longrightarrow
		\mathsf{P}_{2}(\Sigma_b) \longrightarrow \pi_1(\Sigma_b,
		\, p_1)
		\longrightarrow 1.
	\end{equation}
\end{proposition}

\begin{figure}[H]
	\begin{center}
		\includegraphics*[totalheight=2 cm]{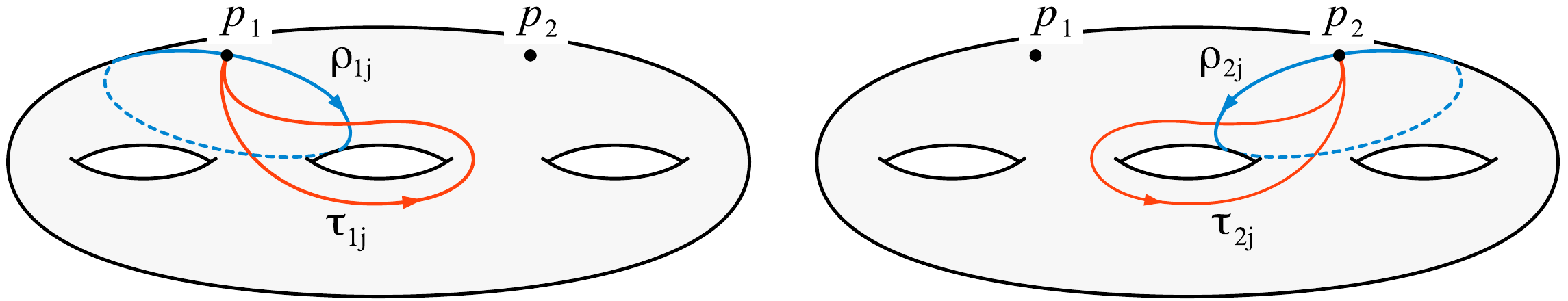}
		\caption{The pure braids $\rho_{1j}, \, \tau_{1j}$,
			$\rho_{2j}, \, \tau_{2j}$
		on $\Sigma_b$} \label{fig1}
	\end{center}
\end{figure}

\medskip
For all $j \in \{1, \ldots, b\}$, let us consider now the $2b$ elements
\begin{equation} \label{eq:braid-generators}
	\rho_{1j}, \; \tau_{1 j}, \; \rho_{2j}, \; \tau_{2 j}
\end{equation}
of $\mathsf{P}_{2}(\Sigma_b)$ represented by the pure braids shown in
Figure \ref{fig1}.

\begin{figure}[H]
	\begin{center}
		\includegraphics*[totalheight=1 cm]{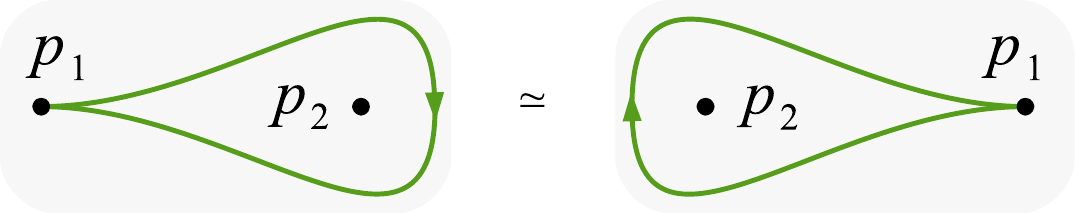}
		\caption{The pure braid $A_{12}$ on $\Sigma_b$} \label{fig2}
	\end{center}
\end{figure}

If $\ell \neq i$, the path corresponding to $\rho_{ij}$ and $\tau_{ij}$
based at $p_{\ell}$ is the constant path. Moreover, let $A_{12}$ be
the pure braid shown in Figure \ref{fig2}. In terms of the isomorphism
\eqref{eq:iso-braids}, the generators $\rho_{ij}$, $\tau_{ij}$ correspond to
the generators of $\pi_1(\Sigma_b \times \Sigma_b - \Delta, \, \mathscr{P})$
coming from the usual description of $\Sigma_b$ as the identification space
of a regular $2b$-gon, whereas $A_{12}$ corresponds to the homotopy class
in $\Sigma_b \times \Sigma_b - \Delta$ of a topological loop in $\Sigma_b
\times \Sigma_b$ that ``winds once'' around $\Delta$.

The elements
\begin{equation} \label{eq:gens-kernel}
	\rho_{21},\ldots, \rho_{2b}, \; \tau_{21},\ldots, \tau_{2b}, \; A_{12}
\end{equation}
can be seen as generators of the kernel $\pi_1(\Sigma_b - \{p_1\}, \, p_2)$
in \eqref{eq:split-braid}, whereas the elements
\begin{equation} \label{eq:gens-quotient}
	\rho_{11},\ldots, \rho_{1b}, \; \tau_{11},\ldots, \tau_{1b}
\end{equation}
are lifts of a set of generators of $\pi_1(\Sigma_b, \, p_1)$ via the quotient
map $\mathsf{P}_{2}(\Sigma_b) \to \pi_1(\Sigma_b, \, p_1)$, namely, they
form a complete system of coset representatives for $\pi_1(\Sigma_b, \, p_1)$.

By Proposition \ref{prop:split-braid}, the braid group
$\mathsf{P}_{2}(\Sigma_b)$ is a semi-direct product of the two groups
$\pi_1(\Sigma_b - \{p_1\}, \, p_2)$ and $\pi_1(\Sigma_b, \, p_1)$,
whose presentations are both well-known; then, in order to write down a
presentation for $\mathsf{P}_{2}(\Sigma_b)$, it only remains to specify
how the generators in \eqref{eq:gens-quotient} act by conjugation on
those in \eqref{eq:gens-kernel}. This is provided by the following result,
cf. \cite[Theorem 1.6]{CaPol19}, where the conjugacy relations are expressed
in the commutator form (i.e., instead of $xyx^{-1}=z$ we write $[x, \,
y]=zy^{-1}$).

\begin{proposition}[{\cite[Theorem 7]{GG04}}] \label{thm:presentation-braid}
	The group $\mathsf{P}_{2}(\Sigma_b)$ admits the following
	presentation. \\ \\
	\emph{Generators}

	$\rho_{1j}, \; \tau_{1j}, \;  \rho_{2j}, \; \tau_{2 j}, \; A_{1
	2} \quad
	j=1,\ldots, b.$ \\\\
	\emph{Relations}
	\begin{itemize}
		\item \emph{Surface relations:}
			\begin{align} \label{eq:presentation-0}
				& [\rho_{1b}^{-1}, \, \tau_{1b}^{-1}] \,
				\tau_{1b}^{-1} \, [\rho_{1
				\,b-1}^{-1}, \, \tau_{1 \,b-1}^{-1}] \,
				\tau_{1\,b-1}^{-1} \cdots
				[\rho_{11}^{-1}, \, \tau_{11}^{-1}] \,
				\tau_{11}^{-1} \, (\tau_{11} \,
				\tau_{12} \cdots \tau_{1b})=A_{12} \\
				& [\rho_{21}^{-1}, \, \tau_{21}] \, \tau_{21}
				\, [\rho_{22}^{-1}, \,
				\tau_{22}] \, \tau_{22}\cdots
				[\rho_{2b}^{-1}, \, \tau_{2b}] \, \tau_{2b}
				\, (\tau_{2b}^{-1} \, \tau_{2 \, b-1}^{-1}
				\cdots \tau_{21}^{-1}) =A_{12}^{-1}
			\end{align}
		\item \emph{Action of} $\rho_{1j}:$
			\begin{align} \label{eq:presentation-1}
				[\rho_{1j}, \, \rho_{2k}]& =1  &  \mathrm{if}
				\; \; j < k \\
				[\rho_{1j}, \, \rho_{2j}]& = 1 & \\
				[\rho_{1j}, \, \rho_{2k}]& =A_{12}^{-1} \,
				\rho_{2k}\, \rho_{2j}^{-1} \,
				A_{12} \, \rho_{2j}\, \rho_{2k}^{-1} \;
				\;&  \mathrm{if} \;  \; j > k \\
				& \\
				[\rho_{1j}, \, \tau_{2k}]& =1 & \mathrm{if}\;
				\; j < k \\
				[\rho_{1j}, \, \tau_{2j}]& = A_{12}^{-1} & \\
				[\rho_{1j}, \, \tau_{2k}]& =[A_{12}^{-1},
				\, \tau_{2k}]  & \mathrm{if}\;
				\; j > k \\
				& \\
				[\rho_{1j}, \,A_{12}]& =[\rho_{2j}^{-1},
				\,A_{12}] &
			\end{align}

		\item \emph{Action of} $\tau_{1j}:$
			\begin{align} \label{eq:presentation-3}
				[\tau_{1j}, \, \rho_{2k}]& =1 & \mathrm{if}\;
				\; j < k \\
				[\tau_{1j}, \, \rho_{2j}]& = \tau_{2j}^{-1}\,
				A_{12}\, \tau_{2j} & \\
				[\tau_{1j}, \, \rho_{2k}]& =[\tau_{2j}^{-1},\,
				A_{12}] \; \; & \mathrm{if}
				\;\; j > k \\
				& \\
				[\tau_{1j}, \, \tau_{2k}]& =1 & \mathrm{if}\;
				\; j < k \\
				[\tau_{1j}, \, \tau_{2j}]& = [\tau_{2j}^{-1},
				\, A_{12}] & \\
				[\tau_{1j}, \, \tau_{2k}]& =\tau_{2j}
				^{-1}\,A_{12}\, \tau_{2j}\,
				A_{12}^{-1}\, \tau_{2k}\,A_{12}\, \tau_{2j}
				^{-1}\,A_{12}^{-1}\,
				\tau_{2j}\,\tau_{2k}^{-1} \; \;  &
				\mathrm{if}\;  \; j > k \\
				&  \\
				[\tau_{1j}, \,A_{12}]& =[\tau_{2j}^{-1},
				\,A_{12}] &
			\end{align}
	\end{itemize}
\end{proposition}

\begin{remark} \label{rmk:no-factorization}
	The inclusion map $\iota \colon \Sigma_b \times \Sigma_b - \Delta
	\to \Sigma_b \times \Sigma_b$ induces a group epimorphism $\iota_*
	\colon \mathsf{P}_2(\Sigma_b) \to \pi_1(\Sigma_b \times \Sigma_b,
	\, \mathscr{P})$, whose kernel is the normal closure of the subgroup
	generated by $A_{12}$. Thus, given any group homomorphism $\varphi
	\colon
	\mathsf{P}_2(\Sigma_b) \to G$, it factors through $\pi_1(\Sigma_b
		\times
	\Sigma_b, \, \mathscr{P})$ if and only if $\varphi(A_{12})$
	is trivial.
\end{remark}

Tedious but straightforward calculations show that the presentation given in
Proposition \ref{thm:presentation-braid} is invariant under the substitutions
\begin{equation*}
	A_{12} \longleftrightarrow A_{12}^{-1}, \quad \tau_{1j}
	\longleftrightarrow
	\tau_{2\; \; b+1-j}^{-1}, \quad  \rho_{1j} \longleftrightarrow
	\rho_{2\;\;
	b+1-j},
\end{equation*}
where $j \in \{1,\ldots, b \}$. These substitutions correspond to the
involution of $\mathsf{P}_2(\Sigma_b)$ induced by a reflection of $\Sigma_b$
switching the $j$-th handle with the $(b+1-j)$-th handle for all $j$. Hence
we can exchange the roles of $p_1$ and $p_2$ in \eqref{eq:split-braid}, and
see $\mathsf{P}_2(\Sigma_b)$ as the middle term of a split short exact
sequence of the form
\begin{equation} \label{eq:split-braid-new}
	1 \longrightarrow \pi_1(\Sigma_b - \{p_2\}, \, p_1)  \longrightarrow
	\mathsf{P}_{2}(\Sigma_b) \longrightarrow \pi_1(\Sigma_b, \, p_2)
	\longrightarrow 1,
\end{equation}
induced by the projection onto the second component
\begin{equation} \label{eq:proj-second}
	(\Sigma_b \times \Sigma_b - \Delta, \, \mathscr{P}) \to (\Sigma_b,
	\, p_2).
\end{equation}
The elements
\begin{equation} \label{eq:gens-kernel-new}
	\rho_{11},\ldots, \rho_{1b}, \; \tau_{11},\ldots, \tau_{1b}, \; A_{12}
\end{equation}
can be seen as generators of the kernel $\pi_1(\Sigma_b - \{p_2\}, \, p_1)$
in \eqref{eq:split-braid-new}, whereas the elements
\begin{equation} \label{eq:gens-quotient-new}
	\rho_{21},\ldots, \rho_{2b}, \; \tau_{21},\ldots, \tau_{2b}
\end{equation}
yield a complete system of coset representatives for $\pi_1(\Sigma_b,
\, p_2)$.

We can now define the objects studied in this paper.

\begin{definition} \label{def:pure-braid-quotient}
	Take positive integers $b, \, n \geq 2$. A finite group $G$ is
	called a
	\emph{pure braid quotient of type} $(b, \, n)$ if there exists a group
	epimomorphism
	\begin{equation} \label{eq:group-epimorphism}
		\varphi \colon \mathsf{P}_2(\Sigma_b) \to G
	\end{equation}
	such that $\varphi(A_{12})$ has order $n$.
\end{definition}

\begin{remark}\label{rmk:purebraidquotient_nonabelian}
	Since we are assuming $n \geq 2$, the element
	$\varphi(A_{12})$ is non-trivial and so the epimorphism $\varphi$ does
	not factor through $\pi_1(\Sigma_b \times \Sigma_b, \, \mathscr{P})$,
	see Remark \ref{rmk:no-factorization}. The geometrical relevance of this
	condition will be explained in Section \ref{sec:ddkf}. The same
	condition also shows that a pure braid quotient is necessarily non-abelian,
	because $\varphi(A_{12})$ is a non-trivial commutator in $G$, see
	\eqref{eq:presentation-1}.
\end{remark}

Sometimes we will use the term \emph{pure braid
quotient} in order to indicate the full datum of the quotient homomorphism
\eqref{eq:group-epimorphism}, instead of the quotient group $G$ alone.

If $G$ is a pure braid quotient, then the two subgroups
\begin{equation} \label{eq:K1-K2}
	\begin{split}
		K_1&:=\langle  \varphi(\rho_{11}), \; \varphi(\tau_{11}),
		\ldots, \varphi(\rho_{1b}), \;
		\varphi(\tau_{1b}), \; \varphi(A_{12})
		\rangle \\
		K_2&:=\langle  \varphi(\rho_{21}), \; \varphi(\tau_{21}),
		\ldots, \varphi(\rho_{2b}), \;
		\varphi(\tau_{2b}), \; \varphi(A_{12})	\rangle
	\end{split}
\end{equation}
are both normal in $G$, and hence there are two short exact sequences
\begin{equation}
	\begin{split}
		&1 \to K_1 \to G \to Q_2 \to 1 \\
		&1 \to K_2 \to G \to Q_1 \to 1,
	\end{split}
\end{equation}
in which the elements $\varphi(\rho_{21}), \; \varphi(\tau_{21}), \ldots,
\varphi(\rho_{2b}), \; \varphi(\tau_{2b})$
yield a complete system of coset representatives for $Q_2$, whereas the
elements
$\varphi(\rho_{11}), \; \varphi(\tau_{11}), \ldots, \varphi(\rho_{1b}),
\; \varphi(\tau_{1b})$ yield a complete system of coset representatives
for $Q_1$.

Let us end this section with the following definition, whose
geometrical meaning will become clear later, see Remark
\ref{rmk:strong-quotient-and-Kodaira-fibrations} of Section \ref{sec:ddkf}.

\begin{definition} \label{def:strong}
	A pure braid quotient $\varphi \colon \mathsf{P}_2(\Sigma_b) \to G$ is
	called \emph{strong} if $K_1=K_2=G$.
\end{definition}

\section{Extra-special groups as pure braid quotients}
\label{sec:extra-special}

We know that $\mathsf{P}_2(\Sigma_b)$ is residually $p$-finite for all prime
number $p \geq 2$, see \cite[pp. 1481-1490]{BarBel09}. This implies that,
for every $p$, we can find a non-abelian finite  $p$-group $G$ that is a pure
braid quotient of type $(b, \, q)$, where $q$ is a power of $p$. However,
it can be tricky to explicitly describe some of these quotients.

In this section we will present a number of results in this direction,
obtained in the series of articles \cite{CaPol19, Pol20, PolSab21}; our
exposition here will closely follow the treatment given in these papers. Let
us start by introducing the following classical definition, see for instance
\cite[p. 183]{Gor07} and \cite[p. 123]{Is08}.

\begin{definition} \label{def:extra-special}
	Let $p$ be a prime number. A finite $p$-group $G$ is called
	\emph{extra-special} if its center $Z(G)$ is cyclic of order $p$
	and the
	quotient $V=G/Z(G)$ is a non-trivial, elementary abelian $p$-group.
\end{definition}

An elementary abelian $p$-group is a finite-dimensional vector space over
the field $\mathbb{Z}_p$, hence it is of the form $V=(\mathbb{Z}_p)^{\dim
V}$ and $G$ fits into a short exact sequence
\begin{equation} \label{eq:extension-extra}
	1 \to \mathbb{Z}_p \to G \to V \to 1.
\end{equation}
Note that, $V$ being abelian, we must have $[G, \, G]=\mathbb{Z}_p$, namely
the commutator subgroup of $G$ coincides with its center. Furthermore,
since the extension \eqref{eq:extension-extra} is central, it cannot be
split, otherwise $G$ would be isomorphic to the direct product of the two
abelian groups $\mathbb{Z}_p$ and $V$, which is impossible because $G$
is non-abelian.  
It can be also proved that, if $G$ is extra-special,
then  $\dim V$ is even and so $|G|=p^{\dim V +1}$ is an odd power of $p$.

For every prime number $p$, there are precisely two isomorphism classes
$M(p)$, $N(p)$ of non-abelian groups of order $p^3$, namely
\begin{equation*}
	\begin{split}
		M(p) & = \langle \mathsf{r}, \, \mathsf{t}, \, \mathsf{z}
		\; | \; \mathsf{r}^p=\mathsf{t}^p=1, \,
		\mathsf{z}^p=1, [  \nr, \,
		\z]=[  \nt,\, \z ]=1, \, [ \nr, \, \nt ]=\z^{-1} \rangle \\
		N(p)& = \langle \nr, \, \nt, \, \z \; | \; \nr^p=\nt^p=\z, \,
		\z^p=1, [ \nr, \,
		\z ]=[ \nt,\, \z ]=1, \, [ \nr, \, \nt ]=\z^{-1} \rangle
	\end{split}
\end{equation*}
and both of them are in fact extra-special, see \cite[Theorem 5.1 of
Chapter 5]{Gor07}.

If $p$ is odd, then the groups $M(p)$ and $N(p)$ are distinguished by
their exponent, which equals $p$ and $p^2$, respectively. If $p=2$, the
group $M(p)$ is isomorphic to the dihedral group $D_8$, whereas $N(p)$
is isomorphic to the quaternion group $\mathsf{Q}_8$.

We can now provide the classification of extra-special $p$-groups, see
\cite[Section 5 of Chapter 5]{Gor07}.
\begin{proposition} \label{prop:extra-special-groups}
	If $b \geq 2$ is a positive integer and $p$ is a prime number, there
	are exactly two isomorphism classes of extra-special $p$-groups
	of order
	$p^{2b+1}$, that can be described as follows.
	\begin{itemize}
		\item The central product $\mathsf{H}_{2b+1}(\mathbb{Z}_p)$
			of $b$ copies
			of $M(p)$, having presentation
			\begin{equation} \label{eq:H5}
				\begin{split}
					\mathsf{H}_{2b+1}(\mathbb{Z}_p)
					= \langle \, & \mathsf{r}_1, \,
					\mathsf{t}_1,
					\ldots, \mathsf{r}_b,\, \mathsf{t}_b,
					\, \z \; |  \; \mathsf{r}_{j}^p =
					\mathsf{t}_{j}^p=\mathsf{z}^p=1, \\
					& [\mathsf{r}_{j}, \, \mathsf{z}]  =
					[\mathsf{t}_{j}, \, \mathsf{z}]= 1, \\
					& [\mathsf{r}_j, \, \mathsf{r}_k]=
					[\mathsf{t}_j, \, \mathsf{t}_k] =
					1, \\
					& [\mathsf{r}_{j}, \,\mathsf{t}_{k}]
					=\mathsf{z}^{- \delta_{jk}} \,
					\rangle.
				\end{split}
			\end{equation}
			If $p$ is odd, this group has exponent $p$.
		\item The central product $\mathsf{G}_{2b+1}(\mathbb{Z}_p)$
			of $b-1$ copies
			of $M(p)$ and one copy of $N(p)$,  having presentation
			\begin{equation} \label{eq:G5}
				\begin{split}
					\mathsf{G}_{2b+1}(\mathbb{Z}_p)
					= \langle \, & \mathsf{r}_1, \,
					\mathsf{t}_1,
					\ldots, \mathsf{r}_b,\, \mathsf{t}_b,
					\, \z \; | \;  \mathsf{r}_{b}^p =
					\mathsf{t}_{b}^p=\mathsf{z}, \\
					&\mathsf{r}_{1}^p = \mathsf{t}_{1}^p=
					\ldots = \mathsf{r}_{b-1}^p =
					\mathsf{t}_{b-1}^p=\mathsf{z}^p=1, \\
					& [\mathsf{r}_{j}, \, \mathsf{z}]  =
					[\mathsf{t}_{j}, \, \mathsf{z}]= 1, \\
					& [\mathsf{r}_j, \, \mathsf{r}_k]=
					[\mathsf{t}_j, \, \mathsf{t}_k] =
					1, \\
					& [\mathsf{r}_{j}, \,\mathsf{t}_{k}]
					=\mathsf{z}^{- \delta_{jk}} \,
					\rangle.
				\end{split}
			\end{equation}
			If $p$ is odd, this group has exponent $p^2$.
	\end{itemize}
\end{proposition}

We are now in a position to state our first two results.

\begin{theorem}[{\cite[Section 4]{CaPol19}, \cite{Pol20}}]
	\label{thm:1}
	If $b \geq 2$ is an integer and $p \geq 5$ is a prime number,
	then both
	extra-special  $p$-groups of order $p^{4b+1}$ are pure braid quotients
	of type $(b, \, p)$. All these quotients are non-strong, in fact $K_1$
	and $K_2$ have index $p^{2b}$ in $G$.
\end{theorem}

\begin{theorem}[{\cite[Section 4]{CaPol19}, \cite{Pol20}}]
	\label{thm:2}
	If $b \geq 2$ is an integer and $p$ is a prime number dividing $b+1$,
	then both extra-special $p$-groups of order $p^{2b+1}$ are strong pure
	braid quotients of type $(b, \, p)$.
\end{theorem}

Theorems \ref{thm:1} and \ref{thm:2} were originally proved by
the first author and A. Causin in \cite{CaPol19}, but only in the case
$G=\mathsf{H}_{4b+1}(\mathbb{Z}_p)$ and $G=\mathsf{H}_{2b+1}(\mathbb{Z}_p)$,
respectively, by using some group-cohomological
results related to the structure of the cohomology algebra $H^*(\Sigma_b
\times \Sigma_b - \Delta, \, \mathbb{Z}_p)$. Let us give here a sketch of the
argument, referring the reader to the aforementioned paper for full details.

Assuming $p \geq 3$, we identified $\mathsf{H}_{4b+1}(\mathbb{Z}_p)$ with
the \emph{symplectic Heisenberg group} $\mathsf{Heis}(V, \, \omega)$, where
\begin{equation}
	V=H_1(\Sigma_b \times \Sigma_b - \Delta, \, \mathbb{Z}_p) \simeq
	H_1(\Sigma_b
	\times \Sigma_b, \, \mathbb{Z}_p) \simeq (\mathbb{Z}_p)^{4b}
\end{equation}
and $\omega$ is a symplectic form on $V$. This group is the central extension
\begin{equation} \label{eq:central}
	1 \to \mathbb{Z}_p \to \mathsf{Heis}(V, \, \omega) \to V \to 1
\end{equation}
of the additive group $V$ given as follows: the underlying set of $
\mathsf{Heis}(V, \, \omega)$ is $V \times \mathbb{Z}_p$, endowed with the
group law
\begin{equation} \label{eq:group-law-Heis}
	(v_1, \, t_1)\,(v_2, \, t_2) = \left(v_1+v_2, \, t_1+t_2 + \frac{1}{2}
	\omega(v_1, \, v_2)\right).
\end{equation}

By basic linear algebra, all symplectic forms on $(\mathbb{Z}_p)^{4b}$ are
equivalent to the standard symplectic form; thus, given two symplectic forms
$\omega_1$, $\omega_2$ on $V$, the two Heisenberg groups $\mathsf{Heis}(V,
\, \omega_1)$, $\mathsf{Heis}(V, \, \omega_2)$ are isomorphic. Moreover,
the center of the Heisenberg group coincides with its commutator subgroup
and is isomorphic to $\mathbb{Z}_p$.

Now, let
\begin{equation}
	\phi \colon \mathsf{P}_2(\Sigma_b) \longrightarrow  V
\end{equation}
be the group epimorphism given by the composition of the reduction mod $p$
map $H_1(\Sigma_b \times \Sigma_b - \Delta, \, \mathbb{Z}) \to V$ with the
abelianization map $\mathsf{P}_2(\Sigma_b) \longrightarrow  H_1(\Sigma_b
\times \Sigma_b - \Delta, \, \mathbb{Z})$. We have a commutative diagram
\begin{equation} \label{diag:lifting-splitting-heisenberg}
%\begin{split}
	\xymatrix{
		&  & & \mathsf{P}_2(\Sigma_b) \ar[d]^{\phi}
		\ar@[red]@{-->}[dl]_{\color{red}{\varphi}} & \\
		1 \ar[r] & \mathbb{Z}_p \ar[r]	& \mathsf{Heis}(V, \,
		\omega) \ar[r]	&
		V
	\ar[r] & 1 }
%\end{split}
\end{equation}
and we denote by $u \in H^2(V, \, \mathbb{Z}_p)$ the cohomology class
corresponding to the bottom Heisenberg extension. Then a lifting $\varphi
\colon \mathsf{P}_2(\Sigma_b) \longrightarrow \mathsf{Heis}(V, \, \omega)$
of $\phi$ exists if and only if  $\phi^*u=0 \in H^2(\mathsf{P}_2(\Sigma_b),
\, \mathbb{Z}_p)$.

The next step is to provide an interpretation of the cohomological condition
$\phi^*u=0$ in terms of the symplectic form $\omega$, and this is achieved
by using the following facts:
\begin{itemize}
	\item we have a natural identification
		\begin{equation} \label{eq:H2}
			H^2(V, \, \mathbb{Z}_p) \simeq \Lambda^2(V^{\vee})
			\oplus V^{\vee}
		\end{equation}
		under which the extension class $u$ giving the Heisenberg
		central extension
		\eqref{eq:central} corresponds to $(\omega, \,
		\epsilon)$. Here $\epsilon
		\colon V \to \mathbb{Z}_p$ stands for the linear functional
		on $V$ defined
		by $\epsilon(v)=w^p$, where $w$ is any preimage of $v$
		in $\mathsf{Heis}(V,
		\, \omega)$;
	\item we have natural identifications
		\begin{equation*}
			V^{\vee} \simeq H^1(\Sigma_b \times \Sigma_b -
			\Delta, \, \mathbb{Z}_p)\simeq
			H^1(\Sigma_b \times \Sigma_b, \, \mathbb{Z}_p)
		\end{equation*}
		and there is a commutative diagram
		\begin{equation} \label{dia:cup-product-1}
			\xymatrix{
				\mathrm{Alt}^2(V) \simeq \wedge^2 V^{\vee}
				\ar@{->}[r]^{\xi}
				\ar@{->}[rd]_{\eta}
				& H^2(\Sigma_b \times \Sigma_b, \,
				\mathbb{Z}_p)  \ar@{->}[d]
				\\
				& H^2(\Sigma_b \times \Sigma_b - \Delta, \,
			\mathbb{Z}_p)}
		\end{equation}
		where the vertical map is the quotient by the $1$-dimensional
		vector
		subspace of $H^2(\Sigma_b \times \Sigma_b, \, \mathbb{Z}_p)$
		generated by
		the class $\delta$ of the diagonal, whereas $\eta$ and $\xi$
		stand for the
		cup product maps;
	\item $\Sigma_b \times \Sigma_b - \Delta$ is an aspherical
		space,
		namely all its higher homotopy group vanish, and so for all $i
		\geq 1$ there is
		a natural isomorphism
		\begin{equation} \label{eq:aspherical}
			H^i(\Sigma_b \times \Sigma_b - \Delta, \,
			\mathbb{Z}_p) \simeq
			H^i(\mathsf{P}_2(\Sigma_b), \, \mathbb{Z}_p)
		\end{equation}
		where $ \mathbb{Z}_p$ is endowed, as an abelian group,
		with the structure
		of trivial $\mathsf{P}_2(\Sigma_b)$-module.
\end{itemize}
Combining all this, we infer that there is a commutative diagram
\begin{equation}
	\begin{split}
		\xymatrix{
			\wedge^2V^{\vee} \oplus V^{\vee} \ar[r]^{\simeq}
			& H^2(V, \, \mathbb{Z}_p)
			\ar[d] \ar[r]^-{\phi^*} & H^2(\mathsf{P}_2(\Sigma_b),
			\, \mathbb{Z}_p)
			\ar[d]^{\simeq} \\
			& \mathrm{Alt}^2(V) \simeq \wedge^2V^{\vee}
			\ar@{->}[r]^-{\eta} & H^2(\Sigma_b
			\times \Sigma_b - \Delta, \, \mathbb{Z}_p),
		}
	\end{split}
\end{equation}
where the isomorphism on the left is \eqref{eq:H2}, the vertical map  on
the left is the projection onto the first summand and the vertical map on
the right is \eqref{eq:aspherical}. Since the projection of the extension
class $u \in H^2(V, \, \mathbb{Z}_p)$ can be naturally identified with
$\omega \in \mathrm{Alt}^2(V)$, we have proved the following

\begin{proposition} \label{prop:interpretation-obstruction}
	The obstruction class $\phi^*u \in H^2(\mathsf{P}_2(\Sigma_b), \,
	\mathbb{Z}_p)$ can be naturally interpreted as the image $\eta(\omega)
	\in
	H^2(\Sigma_b \times \Sigma_b - \Delta, \, \mathbb{Z}_p)$ of the
	symplectic
	form $\omega \in\mathrm{Alt}^2(V)$ via the cup-product map $\eta$.
\end{proposition}

As a consequence, we obtain the following lifting criterion, that we believe
is of independent interest.

\begin{proposition} \label{thm:interpretation-obstruction}
	A lifting $\varphi \colon \mathsf{P}_2(\Sigma_b) \longrightarrow
	\mathsf{Heis}(V, \, \omega)$ of $\phi \colon \mathsf{P}_2(\Sigma_b)
	\to
	V$  exists if and only if $\eta(\omega)=0.$ Furthermore, if $\varphi$
	exists, then $\varphi(A_{12})$ has order $p$ if and only if
	$\xi(\omega)
	\in H^2(\Sigma_b \times \Sigma_b, \, \mathbb{Z}_p)$ is a non-zero
	integer
	multiple of the diagonal class $\delta$. In this case, $\varphi$ is
	necessarily surjective.
\end{proposition}
Inspired by Proposition \ref{thm:interpretation-obstruction}, we say that
a symplectic form $\omega \in\mathrm{Alt}^2(V)$ is \emph{of Heisenberg type}
if $\xi(\omega)$ is a non-zero integer multiple of $\delta$; equivalently,
$\omega$ is of Heisenberg type if $\eta(\omega)=0$ and $\xi(\omega) \neq
0$. By the previous discussion it follows that, if $\omega$ is of
Heisenberg type, $\mathsf{Heis}(V, \, \omega)$
is a pure braid quotient of type $(b, \, p)$.

We are therefore left with the task of constructing symplectic forms of
Heisenberg type on $V$. We denote by $\alpha_1, \, \beta_1, \ldots, \alpha_b,
\,  \beta_b$ the images in $H^1(\Sigma_b, \, \mathbb{Z}_p)=H^1(\Sigma_b, \,
\mathbb{Z}) \otimes \mathbb{Z}_p$ of the elements of a basis of $H^1(\Sigma_b,
\, \mathbb{Z})$ which is symplectic with respect to the cup product; then,
we can choose for $V$ the ordered basis
\begin{equation} \label{eq:ordered-basis-V}
	{r}_{11}, \; {t}_{11}, \ldots, r_{1b}, \; {t}_{1b}, \; \;
	{r}_{21}, \;
	{t}_{21}, \ldots, {r}_{2b}, \; {t}_{2b}
\end{equation}
where, under the isomorphism $V \simeq H_1(\Sigma_b \times \Sigma_b, \,
\mathbb{Z}_p)$ induced by the inclusion $\iota \colon \Sigma_b \times \Sigma_b
- \Delta \to \Sigma_b \times \Sigma_b$, the elements ${r}_{1j}$, ${t}_{1j}$,
${r}_{2j}$, ${t}_{2j} \in V$ are the  duals of the elements $\alpha_j
\otimes 1$, $\beta_j \otimes 1$, $1 \otimes \alpha_j$, $1 \otimes \beta_j
\in H^1(\Sigma_b \times \Sigma_b, \, \mathbb{Z}_p) \simeq H^1(\Sigma_b, \,
\mathbb{Z}_p) \otimes H^1(\Sigma_b, \, \mathbb{Z}_p)$, respectively.

Since $p \geq 5$, we can find non-zero scalars $\lambda_1, \ldots, \lambda_b$,
$\mu_1, \ldots, \mu_b \in \mathbb{Z}_p$ such that $1-\lambda_i\mu_1 \neq 0 $
for all $i \in \{1, \ldots, b \}$ and
\begin{equation} \label{eq:lambda-mu}
	\sum_{j=1}^b \lambda_j = \sum_{j=1}^b \mu_j =1.
\end{equation}
Then we consider the alternating form $\omega \colon V \times V
\to \mathbb{Z}_p$ represented, with respect to the ordered basis
\eqref{eq:ordered-basis-V}, by the skew-symmetric matrix
\begin{equation} \label{eq:omega}
	\Omega_b=
	\begin{pmatrix}
		L_b & J_b \\
		J_b & M_b
	\end{pmatrix} \in \mathrm{Mat}(4b, \, \mathbb{Z}_p)
\end{equation}
where the blocks are the elements of $\mathrm{Mat}(2b, \, \mathbb{Z}_p)$
given by
\begin{equation} \label{eq:L-M}
	L_b=\begin{pmatrix}
		\begin{matrix}0 & \lambda_1 \\ - \lambda_1 & 0\end{matrix}
		& & 0 \\
		& \ddots & \\
		0 & & \begin{matrix}0 & \lambda_b \\ - \lambda_b & 0
		\end{matrix}
	\end{pmatrix} \quad \quad
	M_b=\begin{pmatrix}
		\begin{matrix}0 & \mu_1 \\ - \mu_1 & 0\end{matrix} & & 0 \\
		& \ddots & \\
		0 & & \begin{matrix}0 & \mu_b \\ - \mu_b & 0
		\end{matrix}
	\end{pmatrix}
\end{equation}
\begin{equation} \label{eq:J}
	J_b=\begin{pmatrix}
		\begin{matrix}
		0 & -1 \\ 1 & \; \; 0\end{matrix} & & 0 \\
		& \ddots & \\
		0 & &
		\begin{matrix}0 & -1 \\ 1 & \; \; 0
		\end{matrix}
	\end{pmatrix}
\end{equation}
Standard Gaussian elimination shows that
\begin{equation} \label{eq:determinant-Omega}
	\det \Omega_b = (1-\lambda_1 \mu_1)^2 (1- \lambda_2 \mu_2)^2 \cdots
	(1-\lambda_b \mu_b)^2 >0
\end{equation}
and so $\omega$ is non-degenerate. Moreover, a direct computation
yields $\xi(\omega)= \delta$, that is, $\omega$ is of Heisenberg
type. The calculation of the indices of $K_1$ and $K_2$ in $G$ is now
straightforward, and this completes the proof of Theorem  \ref{thm:1}
in the case $G=\mathsf{H}_{4b+1}(\mathbb{Z}_p)$.

\medskip
Now, let us assume that $p$ divides $b+1$, so that $-b=1$ holds in $\mathbb{Z}_p$,
and take
\begin{equation}
	\lambda_1=\ldots=\lambda_b=\mu_1=\ldots =\mu_b = -1 \in \mathbb{Z}_p.
\end{equation}
Therefore relations \eqref{eq:lambda-mu} are satisfied and the
same computations as in the previous case show that the corresponding
alternating form $\omega$ satisfies $\xi(\omega)=\delta$. However,
$\omega$ is not symplectic, since its associate matrix
\begin{equation} \label{eq:omega-degenerate}
	\Omega_b=
	\begin{pmatrix}
		J_b & J_b \\
		J_b & J_b
	\end{pmatrix} \in \mathrm{Mat}(4b, \, \mathbb{Z}_p)
\end{equation}
has rank $2b$ and, subsequently, $\omega$ has a $2b$-dimensional kernel
$V_0$, namely
\begin{equation} \label{eq:ker-omega}
	V_0 = \langle r_{11}-r_{21}, \, t_{11}-t_{21}, \ldots, r_{1b}-r_{2b},
	\,
	t_{1b}-t_{2b} \rangle.
\end{equation}
The set $V \times \mathbb{Z}_p$, with the operation
\eqref{eq:group-law-Heis}, is a group whose center equals $V_0 \times
\mathbb{Z}_p$ and that, with slight abuse of notation, we denote again by
$\mathsf{Heis}(V, \, \omega)$. Furthermore, the argument in Proposition
\ref{thm:interpretation-obstruction} still applies, providing the existence
of a lifting $\mathsf{P}_2(\Sigma_b) \to \mathsf{Heis}(V, \, \omega)$.
Setting $W=V/V_0$, the alternating form
on $V$	descends to a symplectic form on $W$, that we denote it again
by $\omega$; so $\mathsf{Heis}(W, \, \omega)$ is a genuine Heisenberg
group, endowed with a group epimorphism $\mathsf{Heis}(V, \, \omega) \to
\mathsf{Heis}(W, \, \omega)$. Composing this epimorphism with the
lifting $\mathsf{P}_2(\Sigma_b) \to \mathsf{Heis}(V, \, \omega)$, we obtain
a group epimorphism $\varphi \colon
\mathsf{P}_2(\Sigma_b) \to \mathsf{Heis}(W, \, \omega)$ such that
$\varphi(A_{12})$ is non-trivial and central, hence of order $p$.

Since $W$ is a $\mathbb{Z}_p$-vector space of dimension
$2b$, the group $\mathsf{Heis}(W, \, \omega)$ is isomorphic to
$\mathsf{H}_{2b+1}(\mathbb{Z}_p)$; finally, a simple computation based on the expression \eqref{eq:ker-omega} for $\ker \Omega_b$ yields $K_1=K_2=G$, and this shows Theorem \ref{thm:2} in the case
$G=\mathsf{H}_{2b+1}(\mathbb{Z}_p)$.

The proof of Theorems \ref{thm:1} and \ref{thm:2} in full generality (i.e, for all extra-special groups) was given
in \cite{Pol20}, using a completely algebraic technique that avoided the use
of symplectic geometry and of group cohomology. It is based on the following

\begin{definition} \label{def:ddks}
	Let $G$ be a finite group. A \emph{diagonal double Kodaira structure}
	of type $(b, \, n)$ on $G$ is an ordered set of $4b+1$ generators
	\begin{equation} \label{eq:ddks}
		\mathfrak{S}=(\nr_{11}, \, \nt_{11}, \ldots, \nr_{1b}, \,
			\nt_{1b}, \;
		\nr_{21}, \, \nt_{21}, \ldots, \nr_{2b}, \, \nt_{2b}, \, \z),
	\end{equation}
	with $o(\z)=n$, that are images of the ordered set of generators
	\begin{equation}
		(\rho_{11}, \, \tau_{11}, \ldots, \rho_{1b}, \, \tau_{1b},
			\; \rho_{21},
		\, \tau_{21}, \ldots, \rho_{2b}, \, \tau_{2b}, \, A_{12})
	\end{equation}
	via a pure braid quotient $\varphi \colon \mathsf{P}_2(\Sigma_b)
	\to G$
	of type $(b, \, n)$. The structure  is called \emph{strong} if
	\begin{equation} \label{eq:strong-structure}
		\langle \nr_{11}, \, \nt_{11}, \ldots, \nr_{1b}, \, \nt_{1b}
		\rangle =
		\langle \nr_{21}, \, \nt_{21}, \ldots, \nr_{2b}, \, \nt_{2b}
		\rangle = G.
	\end{equation}
\end{definition}
Therefore, checking whether $G$ is a pure braid quotient of type $(b, \,
n)$ is equivalent to checking whether it admits a diagonal double Kodaira
structure $\mathfrak{S}$ of type $(b, \, n)$. Moreover, by definition,
$\varphi \colon \mathsf{P}_2(\Sigma_b)
\to G$ is strong if and only if $\mathfrak{S}$ is.

Let us refer now to the presentations for extra-special $p$-groups given in Proposition \ref{prop:extra-special-groups}. Assuming that $p$ divides $b+1$, in both cases
$G=\mathsf{H}_{2b+1}(\mathbb{Z}_p)$ and $G=\mathsf{G}_{2b+1}(\mathbb{Z}_p)$
we can obtain a strong diagonal double Kodaira structure $\mathfrak{S}$
on $G$ by setting
\begin{equation*}
	\nr_{1j}=\nr_{2j}=\nr_j, \quad \nt_{1j}=\nt_{2j}=\nt_j
\end{equation*}
for all $j=1, \ldots, b$. The divisibility condition is necessary to ensure
that the element of $\mathfrak{S}$ satisfy the two relations coming  from
the surface relations in Proposition \ref{thm:presentation-braid}. This proves Theorem
\ref{thm:2}.

In order to prove Theorem \ref{thm:1}, it is convenient to consider the
following alternative presentation of extra-special $p$-groups. Consider
any non-degenerate, skew-symmetrix matrix $A=(a_{jk})$ of order $2b$ over
$\mathbb{Z}_p$, and consider the finitely presented groups
\begin{equation} \label{eq:H(A)}
	\begin{split}
		\mathsf{H}(A) = \langle \,  \mathsf{x}_1, \ldots,
		\mathsf{x}_{2b}, \, \mathsf{z} \; | \; & \mathsf{x}_{1}^p
		= \ldots
		=\mathsf{x}_{2b}^p=\mathsf{z}^p=1, \\
		& [\x_1, \, \z] = \ldots = [\x_{2b}, \, \z]=1,\\
		& [\x_j, \, \x_k] =\z^{a_{jk}} \, \rangle,
	\end{split}
\end{equation}

\begin{equation} \label{eq:G(A)}
	\begin{split}
		\mathsf{G}(A) = \langle \,  \mathsf{x}_1, \ldots,
		\mathsf{x}_{2b},
		\, \mathsf{z} \; | \; & \mathsf{x}_{1}^p  = \ldots
		=\mathsf{x}_{2b-2}^p=\mathsf{z}^p=1, \\
		& \x_{2b-1}^p  =\x_{2b}^p= \z, \\
		& [\x_1, \, \z]  = \ldots = [\x_{2b}, \, \z]=1,\\
		& [\x_j, \, \x_k]  =\z^{a_{jk}} \, \rangle,
	\end{split}
\end{equation}
where the exponent in $\z^{a_{jk}}$ stands for any representative
in $\mathbb{Z}$ of $a_{jk} \in \mathbb{Z}_p$. Standard computations
show that $\mathsf{H}(A) \simeq \mathsf{H}_{2b+1}(\mathbb{Z}_p)$
and $\mathsf{G}(A) \simeq \mathsf{G}_{2b+1}(\mathbb{Z}_p)$. Now we
can take as $A$ the matrix $\Omega_b \in \operatorname{Mat}(4b, \,
\mathbb{Z}_p)$ given in \eqref{eq:omega}. Setting $G=\mathsf{H}(\Omega_b)$
or $G=\mathsf{G}(\Omega_b)$, the group $G$ is generated by a set of $4b+1$
elements
\begin{equation}
	\S = \{\nr_{11}, \, \nt_{11}, \ldots, \nr_{1b}, \, \nt_{1b}, \;
		\nr_{21},
	\, \nt_{21}, \ldots, \nr_{2b}, \, \nt_{2b}, \; \z \}
\end{equation}
subject to the relations \eqref{eq:H(A)} or \eqref{eq:G(A)}, respectively. One
can check that $\mathfrak{S}$ provides a diagonal double Kodaira structure
of type $(b, \, p)$ on $G$, and so a diagonal double Kodaira structure of
the same type on the isomorphic group $\mathsf{H}_{4b+1}(\mathbb{Z}_p)$
or $\mathsf{G}_{4b+1}(\mathbb{Z}_p)$. This proves Theorem \ref{thm:1}.

\begin{remark} \label{rmk:thm-1-2-order}
	In particular, the pure braid quotients of smallest order 
	detected by the methods detailed so far are the extra-special 
	groups of order $2^7=128$, 
	corresponding to the case $(b, \, p)=(3, \, 2)$ in 
	Theorem \ref{thm:2}. 
\end{remark}

Recently, in the paper \cite{PolSab21}, we were able to significantly lower the value of $|G|$,
actually providing a sharp lower bound for the order of a pure braid quotient.

\begin{theorem}[{\cite{PolSab21}}] \label{thm:3}
	Assume that $G$ is a finite group that is a pure braid quotient. Then
	$|G|
	\geq 32$, with equality if and only if $G$ is extra-special. In
	this case,
	the following holds.
	\begin{itemize}
		\item There are precisely $2211840=1152 \cdot 1920$
			distinct group
			epimorphisms $\varphi \colon \mathsf{P}_2(\Sigma_2)
			\to G$, and all of them
			make  $G$ a strong pure braid quotient of type $(2,
			\, 2)$.
		\item if $G=G(32, \, 49)=\mathsf{H}_5(\mathbb{Z}_2)$,
			these epimorhisms
			form $1920$ orbits under the natural action of
			$\operatorname{Aut}(G)$.
		\item if $G=G(32, \, 50)= \mathsf{G}_5(\mathbb{Z}_2)$,
			these epimorhisms
			form $1152$ orbits under the natural action of
			$\operatorname{Aut}(G)$.
	\end{itemize}
\end{theorem}

The proof of Theorem \ref{thm:3} is obtained again by looking at the
diagonal double Kodaira structures on $G$, see Definition \ref{def:ddks}. 

\begin{remark} \label{rmk:purebraidquotients_nonCCT}
	A key observation is that if $G$ is
	a CCT-group, namely \( G \) is not abelian and 
	commutativity is a transitive relation on the set
	of the non-central elements, then $G$ admits no 
	diagonal double Kodaira
	structures and, subsequently, it cannot be a pure braid quotient.
\end{remark}

A long but straightforward analysis shows that there are precisely eight
non-CCT groups with $G \leq 32$, namely $G=\mathsf{S}_4$ and $G=G(32, \,
t)$ with $t \in \{6, \, 7, \, 8, \, 43, \, 44, \, 49, \, 50\}$.

These case are handled separately, and a refined analysis proves that
only $G(32, \, 49)$ and $G(32, \, 50)$, i.e. the two extra-special groups,
admit diagonal double Kodaira structures.

Finally, the number of such structures in each case is computed by using the
same techniques as in \cite{Win72}; more precisely, we exploit the fact that
$V=G/Z(G)$ can be endowed with a natural structure of $4$-dimensional symplectic vector space over $\mathbb{Z}_2$,
and that $\operatorname{Out}(G)$ embeds in $\mathrm{Sp}(4, \, \mathbb{Z}_2)$
as the orthogonal group associated with the quadratic form $q$ on $V$ related
to the symplectic form $(\cdot \;, \cdot)$ by $q(\overline{\mathsf{x}}\,
\overline{\mathsf{y}})=q(\overline{\mathsf{x}}) + q(\overline{\mathsf{y}})+
(\overline{\mathsf{x}}, \, \overline{\mathsf{y}})$.

\section{Geometrical application: diagonal double Kodaira fibrations}
\label{sec:ddkf}

Recall that a \emph{Kodaira fibration} is a smooth, connected holomorphic
fibration $f_1 \colon S \to B_1$, where $S$ is a compact complex surface
and $B_1$ is a compact complex curve, which is not isotrivial (this means
that not all fibres are biholomorphic each other). The genus
$b_1:=g(B_1)$ is called the \emph{base genus} of the fibration, whereas
the genus $g:=g(F)$, where $F$ is any fibre, is called the \emph{fibre genus}.
\begin{definition} \label{def:double-kodaira}
	A \emph{double Kodaira surface} is a compact complex surface $S$,
	endowed
	with a \emph{double Kodaira fibration}, namely a surjective,
	holomorphic map
	$f \colon S \to B_1 \times B_2$ yielding, by composition with
	the natural
	projections, two Kodaira fibrations $f_i \colon S \to B_i$, $i=1,
	\,2$.
\end{definition}

With a slight abuse of notation, in the sequel we will use the symbol
$\Sigma_b$ to indicate both a closed Riemann surface of genus $b$ and its
underlying real surface. If a finite group $G$ is a pure braid
quotient of type $(b, \, n)$ then, by using Grauert-Remmert's extension theorem together with Serre's GAGA, the group epimorphism $\varphi \colon
\mathsf{P}_2(\Sigma_b) \to G$ yields the existence of a smooth, complex,
projective surface $S$ endowed with a Galois cover
\begin{equation}
	\mathbf{f} \colon S \to \Sigma_b \times \Sigma_b
\end{equation}
with Galois group $G$ and branched precisely over $\Delta$ with branching
order $n$, see \cite[Proposition 4.4]{CaPol19}. Composing the group
monomorphisms $\pi_1(\Sigma_b - \{p_i\}, \, p_j)  \longrightarrow
\mathsf{P}_{2}(\Sigma_b)$ with $\varphi \colon \mathsf{P}_2(\Sigma_b)
\to G$, we get two homomorphisms
\begin{equation} \label{eq:varphi-i}
	\varphi_1 \colon \pi_1(\Sigma_b -\{p_2\}, \, p_1) \to G, \quad
	\varphi_2
	\colon \pi_1(\Sigma_b -\{p_1\}, \, p_2) \to G,
\end{equation}
whose images are the normal subgroups $K_1$ and $K_2$ defined in $\eqref{eq:K1-K2}$.

By construction, these are the homomorphisms induced by the restrictions
$\mathbf{f}_i \colon \Gamma_i \to \Sigma_b$  of the Galois cover $\mathbf{f}
\colon S \to \Sigma_b \times \Sigma_b$ to the fibres of the two natural
projections $\pi_i
\colon \Sigma_b \times \Sigma_b \to \Sigma_b$. Since $\Delta$ intersects
transversally at a single point all the fibres of the natural projections,
it follows that both such restrictions are branched at precisely one point,
and the number of connected components of the smooth curve
$\Gamma_i \subset S$ equals the index $m_i:=[G : K_i]$ of $K_i$ in $G$.

So, taking the Stein factorizations of the compositions $\pi_i \circ
\mathbf{f} \colon S \to \Sigma_b$ as in the diagram below
\begin{equation} \label{dia:Stein-Kodaira-gi}
	\begin{tikzcd}
		S \ar{rr}{\pi_i \circ \mathbf{f}}  \ar{dr}{f_i} & &
		\Sigma_b    \\
		& \Sigma_{b_i} \ar{ur}{\theta_i} &
	\end{tikzcd}
\end{equation}
we obtain two distinct Kodaira fibrations $f_i \colon S \to \Sigma_{b_i}$,
hence a double Kodaira fibration by considering the product morphism
\begin{equation}
	f=f_1 \times f_2 \colon S \to \Sigma_{b_1} \times \Sigma_{b_2}.
\end{equation}
\begin{definition} \label{def:diagonal-double-kodaira-fibration}
	We call $f \colon S \to \Sigma_{b_1} \times \Sigma_{b_2}$ the
	\emph{diagonal
	double Kodaira fibration} associated with the pure braid quotient
	$\varphi\colon \mathsf{P}(\Sigma_b) \to G$. Conversely, we will
	say that a
	double Kodaira fibration $f \colon S \to \Sigma_{b_1} \times
	\Sigma_{b_2}$
	is \emph{of diagonal type} $(b, \, n)$ if there exists a pure braid
	quotient  $\varphi\colon \mathsf{P}(\Sigma_b) \to G$ of the same
	type such that $f$ is associated with $\varphi$.
\end{definition}

\medskip

Since the morphism $\theta_i \colon \Sigma_{b_i} \to \Sigma_b$ is \'{e}tale
of degree $m_i$, by using the Hurwitz formula we obtain
\begin{equation} \label{eq:expression-gi}
	b_1 -1 =m_1(b-1), \quad  b_2 -1 =m_2(b-1).
\end{equation}
Moreover, the fibre genera $g_1$, $g_2$ of the Kodaira fibrations $f_1
\colon S \to \Sigma_{b_1}$, $f_2 \colon S \to \Sigma_{b_2}$ are computed
by the formulae
\begin{equation} \label{eq:expression-gFi}
	2g_1-2 = \frac{|G|}{m_1} (2b-2 + \mathfrak{n} ), \quad 2g_2-2 =
	\frac{|G|}{m_2} \left( 2b-2 + \mathfrak{n} \right),
\end{equation}
where $\mathfrak{n}:= 1 - 1/n$. Finally, the surface $S$ fits into a diagram
\begin{equation} \label{dia:degree-f-general}
	\begin{tikzcd}
		S \ar{rr}{\mathbf{f}}  \ar{dr}{f} & & \Sigma_b \times
		\Sigma_b  \\
		& \Sigma_{b_1} \times \Sigma_{b_2} \ar[ur, "\theta_1 \times
		\theta_2"{sloped, anchor=south}] &
	\end{tikzcd}
\end{equation}
so that the diagonal double Kodaira fibration $f \colon S \to  \Sigma_{b_1}
\times \Sigma_{b_2}$ is a finite cover of degree $\frac{|G|}{m_1m_2}$,
branched precisely over the curve
\begin{equation} \label{eq:branching-f}
	(\theta_1 \times \theta_2)^{-1}(\Delta)=\Sigma_{b_1} \times_{\Sigma_b}
	\Sigma_{b_2}.
\end{equation}
Such a curve is always smooth, being the preimage of a smooth divisor via an
\'{e}tale morphism. However, it is reducible in general, see \cite[Proposition
4.11]{CaPol19}.

\begin{remark} \label{rmk:strong-quotient-and-Kodaira-fibrations}
	By definition, the pure braid quotient $\varphi \colon
	\mathsf{P}_2(\Sigma_b)
	\to G$ is strong (see Definition \ref{def:strong})
	if and only if $m_1=m_2=1$, that in turn implies $b_1=b_2=b$,
	i.e.,
	$f=\mathbf{f}$. In other words, $\varphi$ is strong if and only
	if no Stein
	factorization as in \eqref{dia:Stein-Kodaira-gi} is needed or,
	equivalently,
	if and only if the Galois cover $\mathbf{f}\colon S \to \Sigma_b
	\times
	\Sigma_b$ induced by $\varphi$ is already a double Kodaira
	fibration,
	branched on the diagonal $\Delta \subset \Sigma_b \times \Sigma_b$.
\end{remark}

We can now compute the invariants of $S$ as follows, see \cite[Proposition
4.8]{CaPol19}.

\begin{proposition} \label{prop:invariant-S-G}
	Let $f \colon S \to \Sigma_{b_1} \times \Sigma_{b_2}$ be a diagonal
	double
	Kodaira fibration, associated with a pure braid quotient $\varphi
	\colon \mathsf{P}_2(\Sigma_b) \to G$
	of type $(b, \, n)$. Then we have
	\begin{equation} \label{eq:invariants-S-G}
		\begin{split}
			c_1^2(S) & = |G|\,(2b-2) ( 4b-4 + 4 \mathfrak{n}
			- \mathfrak{n}^2 ) \\
			c_2(S) & =   |G|\,(2b-2) (2b-2 + \mathfrak{n})
		\end{split}
	\end{equation}
	where $\mathfrak{n}=1-1/n$.	As a consequence, the slope and
	the signature
	of $S$ can be
	expressed as
	\begin{equation} \label{eq:slope-signature-S-G}
		\begin{split}
			\nu(S) & = \frac{c_1^2(S)}{c_2(S)} = 2+ \frac{2
				\mathfrak{n} -
			\mathfrak{n}^2}{2b-2 + \mathfrak{n} } \\
			\sigma(S)& = \frac{1}{3}\left(c_1^2(S) - 2 c_2(S)
			\right)
			=\frac{1}{3}\,|G|\,(2b-2)\left(1-\frac{1}{n^2}\right).
		\end{split}
	\end{equation}
\end{proposition}

\begin{remark} \label{rmk:not-all-diagonal}
	Not all double Kodaira fibrations are of diagonal type. In fact,
	if $S$
	is of diagonal type, then its slope satisfies $\nu(S)=2 +s$, where $
	0 <
	s < 6-4 \sqrt{2}$, see \cite{Pol20}.

\end{remark}

We can now specialize these results, by taking as $G$ an extra-special
$p$-group and using what we have proved in Section \ref{sec:extra-special}.

Fix $b=2$ and let $p \geq 5$ be a prime number. Then every
extra-special $p$-group $G$ of order $p^{4b+1}=p^9$ is a non-strong pure
braid quotient of type $(2, \, p)$ and such that $m_1=m_2=p^{2b}$, see Theorem
\ref{thm:1}. Setting $b':=p^{4}+1$, cf. equations \eqref{eq:expression-gi},
by \cite[Proposition 4.11]{CaPol19} the associated diagonal double Kodaira
fibration $f \colon S \to \Sigma_{b'} \times \Sigma_{b'}$ is a cyclic cover
of degree $p$, branched over a reduced, smooth divisor $D$ of the form
\begin{equation}
	D=\sum_{c \in (\mathbb{Z}_p)^{2b}} D_c
\end{equation}
where the $D_c$ are pairwise disjoint graphs of automorphisms of
$\Sigma_{b'}$.

By using Proposition \ref{prop:invariant-S-G}, we can now construct infinitely
many double Kodaira fibrations with slope strictly higher than $2 + 1/3$.

\begin{theorem} [{\cite[Proposition 4.12]{CaPol19}}, {\cite{Pol20}}] \label{thm:slope}
	Let $f \colon S_p  \to \Sigma_{b'} \times \Sigma_{b'}$ be the diagonal
	double Kodaira fibration associated with a non-strong pure braid
	quotient
	$\varphi \colon \mathsf{P}_2(\Sigma_2) \to G$ of type  $(2, \,
	p)$, where
	$G$ is an extra-special $p$-group $G$ of order $p^9$. Then the maximum
	slope $\nu(S_p)$ is attained for precisely two values of $p$, namely
	\begin{equation}
		\nu(S_{5}) = \nu(S_{7})= 2 + \frac{12}{35}.
	\end{equation}
	Furthermore, $\nu(S_{p}) > 2 + 1/3$ for all $p \geq 5$. More
	precisely,
	if $p \geq 7$ the function $\nu(S_{p})$ is strictly decreasing and
	\begin{equation*}
		\lim_{p \rightarrow +\infty} \nu(S_{p}) = 2 + \frac{1}{3}.
	\end{equation*}
\end{theorem}

\begin{remark} \label{rmk:record-slope}
	The original examples by Atiyah, Hirzebruch and Kodaira have slope
	lying in
	the interval $(2, 2 + 1/3]$, see \cite[p. 221]{BHPV03}. Our
	construction
	provides an infinite family of Kodaira fibred surfaces such that
	$2+1/3 <
	\nu(S) \leq 2 + 12/35$, maintaining at the same time a complete
	control on
	both the base genus and the signature.
	By contrast, the ingenious ``tautological construction" used in \cite{CatRol09}
	yields
	a higher slope than ours, namely $2+ 2/3$, but it involves an
	{\'e}tale
	pullback ``of sufficiently large degree'', that completely loses
	control
	on the other quantities. Note that \cite{LLR17} gives (at least in principle) an effective version of the
pullback construction.
\end{remark}

If $p$ is a prime number dividing $b+1$, by Theorem \ref{thm:2} every
extra-special $p$-group $G$ of order $p^{2b+1}$ is a strong pure braid
quotient of type $(b, \, p)$, and this gives in turn a diagonal double
Kodaira fibration $f \colon S \to \Sigma_b \times \Sigma_b$, see Remark
\ref{rmk:strong-quotient-and-Kodaira-fibrations}. If $\boldsymbol{\upomega}
\colon \mathbb{N} \to \mathbb{N}$ stands for the arithmetic function
counting the number of distinct prime factors of a positive integer, see
\cite[p. 335]{HarWr08}, we obtain
\begin{theorem} [{\cite[Corollary 4.18]{CaPol19}}, \cite{Pol20}] \label{thm:lim-sup}
	Let $\Sigma_b$ be any closed Riemann surface of genus $b$. Then
	there exists a
	double Kodaira fibration $f \colon S \to \Sigma_b \times
	\Sigma_b$. Moreover,
	denoting by $\kappa(b)$ the number of such fibrations, we have
	\begin{equation*}
		\kappa(b) \geq \boldsymbol{\upomega}(b+1).
	\end{equation*}
	In particular,
	\begin{equation*}
		\limsup_{b \rightarrow + \infty} \kappa(b) = + \infty.
	\end{equation*}
\end{theorem}
As far as we know, this is the first construction showing that \emph{all}
curves of genus $b\geq 2$ (and not only some special curves with extra
automorphisms) appear in the base of at least one double Kodaira fibration
$f \colon S \to \Sigma_b \times \Sigma_b$. In addition, two Kodaira
fibred surfaces corresponding to two distinct prime divisors of $b+1$ are
non-homeomorphic, because the corresponding signatures are different: just
use \eqref{eq:slope-signature-S-G} with $n=p$ and $|G|=p^{2b+1}$ and note
that, for fixed $b$,
the function expressing $\sigma(S)$ is strictly increasing in $p$. This
shows that the number of topological types of $S$, for a fixed base $\Sigma_b$, can be arbitrarily large.

\medskip
Let us now consider Theorem \ref{thm:3}, whose geometrical translation is
\begin{theorem}[\cite{PolSab21}] \label{thm:3-geometric}
	Let $G$ be a finite group and $\mathbf{f} \colon S \to \Sigma_b \times
	\Sigma_b$ be a Galois cover, with Galois group $G$, branched over
	the diagonal
	$\Delta$ with branching order $n$. Then $|G|\geq 32$, and equality
	holds
	if and only if $G$ is extra-special. In this case, the following
	holds.
	\begin{itemize}
		\item[$\boldsymbol{(1)}$]$\,$ There exist $2211840=1152
			\cdot 1920$ distinct
			$G$-covers $\mathbf{f} \colon S \to \Sigma_2 \times
			\Sigma_2$, and all of
			them are diagonal double Kodaira fibrations such that
			\begin{equation}
				b_1=b_2=2, \quad g_1=g_2=41,  \quad
				\sigma(S)=16.
			\end{equation}
		\item[$\boldsymbol{(2)}$]$\,$ If $G=G(32, \,
			49)=\mathsf{H}_5(\mathbb{Z}_2)$,
			these $G$-covers form
			$1920$ equivalence classes up to cover isomorphisms.
		\item[$\boldsymbol{(3)}$]$\,$ If $G=G(32, \,
			50)=\mathsf{H}_5(\mathbb{Z}_2)$,
			these $G$-covers form
			$1152$ equivalence classes up to cover isomorphisms.
	\end{itemize}
\end{theorem}

As a consequence, we obtain a sharp lower  bound for the signature of a
diagonal double Kodaira fibration. In fact, the second equality in
\eqref{eq:slope-signature-S-G} together with Theorem \ref{thm:3-geometric}
imply that, for every such fibration, we have
\begin{equation}
	\sigma(S) =\frac{1}{3}\,|G|\,(2b-2)\left(1-\frac{1}{\;
	n^2}\right) \\
	\geq \frac{1}{3} \cdot 32 \cdot (2\cdot 2 -2)
	\left(1-\frac{1}{\; 2^2}\right)
	= 16,
\end{equation}
and this in turn establishes the following result.

\begin{theorem}[{\cite{PolSab21}}] \label{thm:signature}
	Let $S$ be a double Kodaira surface, associated with a pure braid quotient $\varphi
	\colon \mathsf{P}_2(\Sigma_b) \to G$ of type $(b, \, n)$. Then
	$\sigma(S) \geq 16$, and equality
	holds precisely when $(b, \, n)=(2, \, 2)$ and $G$ is an extra-special
	group of order $32$.
\end{theorem}

\begin{remark}
If $S$ is a double Kodaira fibration, corresponding to a pure braid quotient $\varphi
	\colon \mathsf{P}_2(\Sigma_b) \to G$ of type $(b, \, n)$, then, using the terminology in \cite{CatRol09}, it is \emph{very simple}. Let us denote by $\mathfrak{M}_S$ the connected component of the Gieseker moduli space of surfaces of general type containing the class of $S$, and by $\mathcal{M}_b$ the moduli space of smooth curves of genus $b$. Thus, by applying \cite[Thm. 1.7]{Rol10} and using the fact that $\Delta \subset \Sigma_b \times \Sigma_b$ is the graph of the identity $\operatorname{id} \colon \Sigma_b \to \Sigma_b$, we infer that every surface in  $\mathfrak{M}_S$ is still a very simple double Kodaira fibration and that there is a natural map of schemes 
	\begin{equation}
	\mathcal{M}_b \to \mathfrak{M}_S,
	\end{equation} 
which is an isomorphism on geometric points. Roughly speaking,  since the branch locus $\Delta \subset \Sigma_b \times \Sigma_b$ is rigid, all the deformations  of $S$ are realized by deformations of $\Sigma_b \times \Sigma_b$ preserving the diagonal, hence by deformations of $\Sigma_b$, cf. \cite[Proposition 4.22]{CaPol19}. In particular, this shows that  $\mathfrak{M}_S$ is a connected and irreducible component of the Gieseker moduli space.
\end{remark}

Every Kodaira fibred surface $S$ has the structure of a real surface
bundle over a real surface, and so $\sigma(S)$ is divisible
by $4$, see \cite{Mey73}.  If, in addition, $S$ has a spin structure,
i.e. its canonical class is $2$-divisible in $\operatorname{Pic}(S)$,
then $\sigma(S)$ is a positive multiple of $16$ by Rokhlin's theorem,
and examples with $\sigma(S)=16$ are constructed in \cite{LLR17}. It is
not known if there exists a Kodaira fibred surface with $\sigma(S) \leq 12$.

Constructing (double) Kodaira fibrations with small signature is
a rather
difficult problem. As far as we know, before the present work the only
examples
with signature $16$ were the ones listed in \cite[Table 3, Cases
	6.2, 6.6,
6.7 (Type 1), 6.9]{LLR17}. The examples in Theorem \ref{thm:signature}
are new, since both the base genera and the fibre genera are
different from the ones in the aforementioned cases. Our results also show
that \emph{every} curve
of genus $2$ is the base
of a double Kodaira fibration with signature $16$. Thus,  we obtain
two families
of dimension $3$ of such fibrations that, to our
knowledge,
provides the first examples of positive-dimensional families of double
Kodaira
fibrations with small signature.

Theorem \ref{thm:signature} also provide new ``double solutions'' to
a problem, posed by G. Mess and included  in Kirby's problem list in
low-dimensional topology, see
\cite[Problem 2.18 A]{Kir97}, asking what is the smallest number $b$ for
which there exists a real surface bundle over a real surface with base
genus $b$ and non-zero signature. We actually have $b=2$, also for double
Kodaira fibrations.

\begin{theorem}[{\cite{PolSab21}}] \label{thm:real-manifold}
	Let  $S$ be a double
	Kodaira surface, associated with a pure braid quotient $\varphi \colon
	\mathsf{P}_2(\Sigma_b) \to G$ of type $(2, \, 2)$, where $G$ is  an
	extra-special group of order $32$. Then the real manifold $X$ underlying
	$S$ is a closed, orientable $4$-manifold of signature $16$ that
	can be realized as a real surface bundle over a real surface of
	genus $2$, with fibre  genus $41$, in two different ways.
\end{theorem}

It is an interesting question whether $16$ and $41$ are the minimum
possible values
for the signature and the fibre genus of a (non necessary diagonal) double
Kodaira fibration $f \colon
S \to \Sigma_2 \times \Sigma_2$; however, this topic exceeds the scope of
this paper.

\section{Beyond \texorpdfstring{\( | G | = 32 \)}{|G|=32}} \label{sec:new}

This last section contains the new result of this article. As we already
observed, so far we have detailed a rather clear picture regarding pure
braid quotient groups \( G\) and the relative diagonal double 
Kodaira fibrations for \( |G| \le 32 \). Indeed, there is no pure braid
quotient of order strictly less than \( 32\) and for \( |G| = 32 \) we have only the two extra-special groups, see Theorem \ref{thm:3}. Furthermore, Theorem
\ref{thm:2}  provides examples of pure braid quotients starting 
with order equal to \( 128\), see Remark \ref{rmk:thm-1-2-order}, and they are
extra-special groups, too. 
It seems then natural to investigate this matter further for \( |G| > 32
\), for instance, in order to look for non extra-special examples.
In this direction we have obtained the following result that
highlights the existence of a gap between orders \( 32 \) and \( 64 \). 
\begin{theorem} \label{thm:new}
	If $G$ is a finite group with $32 < |G| < 64$, then $G$
	is not a pure braid quotient.
\end{theorem}
Here we just give a sketch of the proof, while the full details will appear
elsewhere. We know that the group $G$ cannnot be abelian, see Remark
\ref{rmk:purebraidquotient_nonabelian}, and we also mentioned that
CCT-groups cannot be pure braid quotients, see Remark
\ref{rmk:purebraidquotients_nonCCT}. 
On the other hand, by \cite{PolSab21} we know
that, if $G$ is a pure braid quotient and admits no proper quotients that
are pure braid quotients,
then $G$ must be monolithic, i.e., the intersection $\operatorname{soc}(G)$
of its non-trivial normal subgroups is non-trivial. In fact,
consider the 
epimorphism $\varphi \colon \mathsf{P}_2(\Sigma_b) \to G$ and assume that
there is a non-trivial, normal subgroup $N$ of $G$ such that $\varphi(A_{12})
\notin N$. Then, composing the projection $G\to G/N$ with $\varphi$, we
obtain a pure braid quotient $\bar{\varphi} \colon \mathsf{P}_2(\Sigma_b)
\to G/N$, which leads to a contradiction. It follows that
$\varphi(A_{12}) \in \operatorname{soc}(G)$, in particular $G$ is monolithic.

By Theorem \ref{thm:3} this implies that, if a pure braid quotient $G$ 
satisfies our assumptions on the order, then $G$ must be monolithic.
A straightforward computer calculation with \verb|GAP4| now shows that
there are precisely two non-abelian groups $G$ with
$32 < |G| < 64$  that are
both non-CCT and monolithic, namely $G(54, \, 5)$ and $G(54, \, 6)$; by
the remarks above, they are the only possible candidates to be pure braid 
quotients in that range for $|G|$.
Finally, a brute force check (again by using \verb|GAP4|) shows that these
groups admit no diagonal double Kodaira structure, proving our assertion.

\end{document}